\documentclass[11pt]{amsart}
\usepackage{amsmath}
\usepackage{latexsym}
\usepackage{amssymb}
\usepackage{amscd}
\input epsf

\usepackage{amsfonts}
\usepackage[all]{xy}

\DeclareFontShape{OT1}{cmr}{b}{n}{<12> cmr12}{}

\newtheorem{proposition}{Proposition}[section]

\newtheorem{lemma}[proposition]{Lemma}
\newtheorem{definition}[proposition]{Definition}
\newtheorem{theorem}[proposition]{Theorem}

\newtheorem{remark}[proposition]{Remark}

\newtheorem{lemma-definition}[proposition]{Lemma-Definition}

\newenvironment{dok}{\par\vspace{-5pt}%
\par\noindent\begingroup%
\leftskip=0em\hspace{0em}{\bf Proof.}}%
{\endgroup\hfill$\Box$}

\newcounter{SNO}
\newcounter{THNO}[section]
\newcounter{tmp}

\def\db#1{ \bD^b({#1})}

\def\perf#1{{\mathfrak P}{\mathfrak e}{\mathfrak r}{\mathfrak f}({#1})}
\def\dsing#1{ \bD_{Sg}({#1})}

\def\dg{DG}
\def\dt{H^0\dg}

\def\abe{\operatorname{Pair}}

\def\Hom{{\mathrm H}{\mathrm o}{\mathrm m}}
\def\Hhom{{\mathbb H}{\mathrm o}{\mathrm m}}
\def\Cok{\operatorname{Cok}}

\def\h#1,#2{{\operatorname{Hom}}({#1}\:,\; {#2})}
\def\Ho#1,#2,#3,#4{{\operatorname{Hom}}^{#1}_{#2}({#3}\:,\; {#4})}
\def\Ex#1,#2,#3,#4{{\operatorname{Ext}}^{#1}_{#2}({#3}\:,\; {#4})}

\def\coh{\operatorname{coh}}
\def\Qcoh{\operatorname{Qcoh}}

\def\lto{\longrightarrow}

\def\D{{\mathcal D}}
\def\F{{\mathcal F}}
\def\E{{\mathcal E}}

\def\H{{\mathcal H}}

\def\O{{\mathcal O}}

\def\T{{\mathcal T}}

\def\P{{\mathcal{P}}}

\def\ZZ{{\mathbb Z}}

\def\bE{{\mathbb E}}
\def\bF{{\mathbb F}}
\def\bG{{\mathbb G}}
\def\bT{{\mathbb T}}
\def\bU{{\mathbb U}}
\def\bL{{\mathbb L}}
\def\bP{{\mathbb P}}
\def\Cone{{\ul{\mathbb{C}{\mathrm{one}}}}}
\def\b1{{\mathbf 1}}

\def\MF{\mathrm{MF}}
\def\Ac{\mathit{Ac}}

\def\bD{{\mathbf D}}
\def\bR{{\mathbf R}}

\def\ds#1{ {#1}_{\rm Sg}}
\def\dhf#1{ {#1}_{\rm hf}}

\def\AA{{\mathbb A}}

\def\ZZ{{\mathbb Z}}

\def\Hom{\operatorname{Hom}}

\def\Ext{\operatorname{Ext}}

\def\ext{\underline{{\mathcal E}xt}}
\def\rhom{\bR\underline{{\mathcal H}om}}

\def\Sing{\operatorname{Sing}}

\def\Coker{\operatorname{Coker}\,}

\def\rH{{\operatorname{H}}}

\def\rZ{{\operatorname{Z}}}

\def\id{{\operatorname{id}}}

\def\ul{\underline}
\def\kk{{\mathbf k}}

\def\mod{\operatorname{mod}\!}

\def\wX{\mathfrak{X}}

\title[]{Matrix factorizations for nonaffine LG--models}

\author[]{Dmitri Orlov}

\address{ Algebra Section, Steklov Mathematical Institute RAN,
Gubkin str. 8, Moscow 119991, RUSSIA}

\email{orlov@mi.ras.ru}

\thanks{
This work was partially supported by  RFBR grants 10-01-93113, 11-01-00336, 11-01-00568, NSh grant 4713.2010.1,
by AG Laboratory HSE, RF government grant, ag. 11.G34.31.0023, and by Simons Center for Geometry and Physics.}

\date{}
\dedicatory{
}

\begin{document}

\begin{abstract}
We propose a natural definition of a category of matrix factorizations
for nonaffine Landau-Ginzburg models. For any LG-model
we construct a fully faithful functor from the category of matrix
factorizations defined in this way to the triangulated category of singularities of the corresponding fiber.
We also show that this functor is an equivalence if
the total space of the LG-model is smooth.
\end{abstract}
\maketitle

\vspace*{-0.2cm}
\section*{Introduction}

In the paper \cite{Tr} we established a connection between categories
of D-branes of type B (B-branes) in affine Landau-Ginzburg models and triangulated
categories of singularities of the singular fibers.  A mathematical
definition for the category of B-branes in affine Landau-Ginzburg
models was proposed by M.~Kontsevich. According to his proposal the
superpotential $W$ will deform complexes of coherent sheaves to ``$W$\!-twisted'' complexes. These are 2-periodic chains of maps of
vector bundles in which the composition of two consequtive maps is no
longer required to be zero, but instead is equal to multiplication by
$W.$ Such chains are called matrix factorizations of $W$ and so
Kontsevich predicted that the category of B-branes should be the
category of matrix factorizations.  Kapustin and Li verified
\cite{KL} the equivalence of this definition with the physics notion
of B-branes in LG-models in the case of the usual quadratic
superpotential and gave physics arguments supporting Kontsevich's
proposal in the case of  general superpotentials.

The triangulated category of singularities is defined for any
noetherian scheme as a quotient of the bounded derived category of
coherent sheaves by the subcategory of perfect complexes.  On the
other hand, the definition of the category of matrix factorization was
given only for affine LG-models. Under the conditions that the total
space $X$ of an LG-model is affine and smooth we showed in the paper
\cite{Tr} that the categories of matrix factorizations are equivalent
to the triangulated categories of singularities of the corresponding
fibers.  In particular, this implies that, in spite of the
fact that the category of B-branes is defined using the total space
$X,$ it depends only on the singular fibers of the superpotential.

In this paper we propose a natural generalization of the definition of
the category of matrix factorizations which makes sense for arbitrary,
not necessarily affine LG-models.  We show that with our general
definition the same result as in affine case holds true as long as the
total space of an LG-model is smooth.  In the case of a singular total
space we have only a full embedding of the category of matrix
factorizations in the triangulated category of the corresponding
fiber.  Studying the difference between these categories in the case
of singular total spaces is certainly very interesting and we plan to
return to this question in the future.

When this paper was finished I was informed that a category of matrix
factorizations in the nonaffine case is also discussed in a forthcoming
paper of Kevin Lin and Daniel Pomerleano \cite{LP}.

I am grateful to Ludmil Katzarkov, Tony Pantev, and Leonid Positselski
for useful discussions and their interest in this work.  I would like to
thank the Simons Center for Geometry and Physics for hospitality and
very stimulating atmosphere during the preparation of this paper.

\section{Triangulated categories of singularities}

Let $X$ be a noetherian scheme over a field $\kk.$ Denote by
$\bD(\Qcoh X)$ the unbounded derived category of quasi-coherent
sheaves on $X.$ The derived category $\bD(\Qcoh X)$ has arbitrary
coproducts.  We can consider the full triangulated subcategory of
compact objects, i.e. such objects $C\in\bD(\Qcoh X)$ for which the
functor $\Hom(C,-)$ commutes with all coproducts. It is known (see
\cite{Ne, BVdB}) that the subcategory of compact objects coincides
with the subcategory of perfect complexes $\perf{X}.$ Recall that a
complex on a scheme is called perfect if it is locally
quasi-isomorphic to a bounded complex of locally free sheaves of
finite type.

For a noetherian scheme it is also natural to consider the abelian
subcategory $\coh(X)\subset\Qcoh(X)$ of coherent sheaves and the
bounded derived categories $\db{\coh X}$ of coherent sheaves on $X.$
The natural functor from $\db{\coh X}$ to the unbounded derived
category of quasi-coherent sheaves $\bD(\Qcoh X)$ is fully faithful
and identifies $\db{\coh X}$ with the full subcategory
$\bD^{\emptyset,\, b}_{\coh}(\Qcoh X)$ consisting of all
cohomologically bounded complexes with coherent cohomology (\cite{Il},
Ex.II, 2.2.2).  In particular $\perf{X}$ can be viewed as a full
triangulated subcategory of $\db{\coh X}.$ For a regular scheme the
inclusion $\perf{X}\subseteq\db{\coh X}$ is an equivalence while for a
singular scheme this inclusion is strict.  The triangulated category
of singularities of $X$ is defined as the Verdier quotient of these
two natural categories that we can attach to a noetherian scheme $X$:

\begin{definition}\label{trcsin} The triangulated category of
  singularities of $X,$ denoted by $\dsing{X},$ is the Verdier
quotient of the bounded derived category $\db{\coh X}$ by the full
triangulated subcategory of perfect complexes $\perf{X}.$
\end{definition}

We say that $X$ satisfies the condition {\sf (ELF)} if
\vspace{6pt}

\begin{itemize}
\item[{\sf (ELF)}] $X$ is separated, noetherian, of finite Krull
  dimension, and has enough locally free sheaves.
\end{itemize}
\

\vspace{1pt}

The last condition means that for any coherent sheaf $\F$ there is an
epimorphism $\E\twoheadrightarrow\F$ with a locally free sheaf $\E.$
It also implies that any perfect complex is globally (not only
locally) quasi-isomorphic to a bounded complex of locally free sheaves
of finite type.  For example, any quasi-projective scheme satisfies
these conditions.  Note that if $X$ satisfies (ELF), then any
closed and any open subscheme of $X$ is also noetherian, finite
dimensional and has enough locally free sheaves.  This is automatic
for a closed subscheme, while for an open subscheme $U \subset X$ it
follows from the fact that any coherent sheaf on $U$ can be obtained
as the restriction of a coherent sheaf on $X.$ From now on we will
assume that all schemes we work with satisfy the condition
(ELF).

The definition of triangulated category of singularities extends
verbatim to the case of stacks and noncommutative spaces.  One can
also define a graded version of this category, which is a central
object in the so called Landau-Ginzburg/Calabi-Yau correspondence
discussed in the paper \cite{Grad}.

It is possible to reformulate Definition \ref{trcsin} so that it
will make sense intrinsically in any triangulated category $\D.$ We
say that an object $A\in\D$ is {\sf homologically finite} if for any
object $B\in\D$ all $\Hom(A, B[i])$ are trivial except for finite
number of $i\in \ZZ.$ Such objects form a full triangulated
subcategory  $\dhf{\D} \subset \D.$  We define a triangulated
category $\ds{\D}$ as the quotient $\D/\dhf{\D}.$

It is proved in \cite{Equi} that if a scheme $X$ satisfies (ELF). Then
the subcategory $\dhf{\D}$ of homologically finite objects in
$\D=\db{\coh(X)}$ coincides with the subcategory of perfect complexes
$\perf{X}$ and, hence, $\ds{\D}\cong\dsing{X}.$ Thus the triangulated
category of singularities can be defined internally starting only with
the bounded derived category of coherent sheaves.

A simple but fundamental property of triangulated categories of
singularities is the fact that they are local in the Zariski
topology. Explicitly this means that for any Zariski open embedding $j:
U\hookrightarrow X,$ for which $\Sing(X)\subset U,$ the functor
$\bar{j}^*:\dsing{X}\to \dsing{U}$ is an equivalence of triangulated
categories \cite{Tr}.

On the other hand, two analytically isomorphic singularities can have
non-equivalent (but closely related) triangulated categories of
singularities.  The main reason why such categories can fail to be
equivalent is that a triangulated category of singularities is not
necessarily idempotent complete.

However, for any triangulated category $\D$ we can consider its so
called idempotent completion (or Karoubian envelope) $\overline{\D}.$
This is a category that consists of all kernels of all projectors. It
has a natural structure of a triangulated category and the canonical
functor $\D\to\overline{\D}$ is an exact full embedding \cite{BS}.  In
\cite{Form} we show that for any two schemes $X$ and $X'$ satisfying
(ELF), whose formal completions $\wX$ and $\wX'$ along the
singularities are isomorphic, we have that the idempotent completions
of the triangulated categories of singularities $\overline{\dsing{X}}$
and $\overline{\dsing{X'}}$ are equivalent.

Triangulated categories of singularities have additional good
properties in the case of Gorenstein scheme.  If $X$ is Gorenstein and
has finite dimension, then $\O_X$ is a dualizing complex for $X,$ in
the sense that it has a finite injective dimension as a quasi-coherent
sheaf and the natural map
$$
\F\lto \rhom^{\cdot}(\rhom^{\cdot}(\F, \O_X), \O_X)
$$
is an isomorphism for any coherent sheaf $\F.$
In particular, there is
an integer $n_0$ such that
$\ext^{i}(\F, \O_X)=0$ for each quasi-coherent sheaf
$\F$ and all $i>n_0.$

The following statements and their proofs can be found in \cite{Tr}.
\begin{lemma-definition}\label{rres}{\rm (\cite{Tr}, Lemma 1.19)}
Let $X$ be a Gorenstein scheme satisfying  (ELF). We say that
a coherent sheaf $\F$ is Cohen-Macaulay if the following  equivalent
conditions hold.
\begin{itemize}
\item[1)] The sheaves $\ext^{i}(\F, \O_X)$ are trivial for all
$i>0.$
\item[2)] There is a right locally free resolution
$
0\lto\F\lto\{Q^0\lto Q^1 \lto Q^3\lto \cdots \}.
$
\end{itemize}
\end{lemma-definition}

\begin{proposition}\label{fint}{\rm (\cite{Tr}, Prop. 1.23)}
Let $X$ be a Gorenstein scheme satisfying (ELF).
Then any object $A\in\dsing{X}$ is isomorphic
to the image of a Cohen-Macaulay sheaf.
\end{proposition}

For a Gorenstein scheme $X$ that satisfy condition (ELF) we can
calculate morphisms between objects in the triangulated category of
singularities in terms of morphisms between coherent sheaves on $X$
(see \cite{Tr}, Prop 1.21). In particular, it follows that if the
closed subset $\Sing(X)$ is complete then the spaces of morphisms
between any two objects in the triangulated category of singularities
$\dsing{X}$ are finite dimensional.

\section{Matrix factorizations for non-affine Landau-Ginzburg models}
\label{three}

By a Landau-Ginzburg model we will mean the following data: a scheme
$X$ over a field $\kk$ and a regular function $W$ on $X$ such that the
morphism $W: X\to \AA^1_{\kk}$ is flat. (It is equivalent to say
that the map of algebras
$\kk[x]\to \varGamma(\O_X)$ is an injection.) In general, the data of a LG model
should also include a K\"{a}hler form on $X$ but we will ignore this
extra piece of information since it is not needed for the definition
of B-branes. As it was mentioned before in this paper we assume that
the scheme $X$ satisfies the condition (ELF).

\begin{remark}{\rm
Usually in the definition of a LG-model we ask that $X$ be regular and
our final equivalence result holds precisely under this condition.
However for all other considerations regularity of $X$ is not
necessary. Moreover, it seems that it is very interesting to consider
the case of a singular $X$ as well.}
\end{remark}

With any $\kk$\!-point $w_0\in \AA^1$ we can associate a differential
$\ZZ/2\ZZ$\!-graded category $\dg_{w_0}(X, W),$ an exact category
$\abe_{w_0}(X, W),$ and a triangulated category $\dt_{w_0}(X, W)$ that
is the homotopy category for DG category $\dg_{w_0}(X, W).$

Objects of all these categories are  ordered pairs
$$
\ul{\bE}:=\Bigl(
\xymatrix{
\bE_1 \ar@<0.6ex>[r]^{e_1} &\bE_0 \ar@<0.6ex>[l]^{e_0}
}
\Bigl),
$$ where $\bE_0, \bE_1$ are locally free sheaves of finite type on $X$
and the compositions ${e_0 e_1}$ and ${e_1 e_0}$ are the
multiplications by the element $(W-w_0\cdot \b1)\in \varGamma (\O_X).$

Morphisms from $\ul{\bE}$ to $\ul{\bF}$
in the category $\dg_{w_0}(W)$ form $\ZZ/2\ZZ$\!-graded
complex
$$
{\Hhom}(\ul{\bE},\; \ul{\bF})=\bigoplus_{0\le i,j\le 1} \Hom(\bE_i, \bF_j)
$$
with a natural grading $(i-j)\mod\; 2,$ and with a
differential $D$ acting on  a homogeneous element $p$
of degree $k$ as
$$
D p=f \cdot p-(-1)^{k}p\cdot e .
$$

The space of morphisms $\Hom(\ul{\bE},\; \ul{\bF})$ in the category
$\abe_{w_0}(X, W)$ is the space of morphisms in $\dg_{w_0}(X, W)$
which are homogeneous of degree 0 and commute with the differential.

The space of morphisms in the category $\dt_{w_0}(W)$ is the space of
morphisms in $\abe_{w_0}(X, W)$ modulo null-homotopic morphisms, i.e.
$$
\Hom_{\abe_{w_0}(X, W)}(\ul{\bE},\; \ul{\bF})= \rZ^0(\Hhom(\ul{\bE},\;
\ul{\bF})),
\qquad
\Hom_{\dt_{w_0}(X,W)}(\ul{\bE},\; \ul{\bF})= \rH^0(\Hhom(\ul{\bE},\;
\ul{\bF})).
$$
Thus, a morphism $p:\ul{\bE}\to\ul{\bF}$ in the category
$\abe_{w_0}(X,W)$ is a pair of morphisms
$p_1: \bE_1\to \bF_1$ and $p_0: \bE_0\to \bF_0$ such that
$p_1 e_0=f_0 p_0$ and $f_1 p_1=p_0 e_1.$
The morphism $p$ is null-homotopic if there are two morphisms
$s_0: \bE_0\to \bF_1$ and $s_1:\bE_1\to \bF_0$ such that
$p_1=f_0 s_1 + s_0 e_1$ and $p_0=s_1 e_0 + f_1 s_0.$

It is clear that the category $\abe_{w_0}(X, W)$
is an exact category with respect to componentwise
monomorphisms and epimorphisms.

The category $\dt_{w_0}(X, W)$ can be
endowed with a natural structure of a triangulated category.
To specify it we have to define a translation functor $[1]$
and a class of exact triangles.

The translation functor can be defined as a functor
that takes an object $\ul{\bE}$ to the object
$$
\ul{\bE} [1]=\xymatrix{
\Bigl(
\bE_0 \ar@<0.6ex>[r]^{-e_0} &\bE_1 \ar@<0.6ex>[l]^{-e_1}
}
\Bigl),
$$
i.e. it changes the order of the modules and the signs of the maps,
and takes a morphism $p=(p_0, p_1)$ to the morphism
$p[1]=(p_1, p_0).$
We see that the functor $[2]$ is the identity functor.

For any morphism $p:\ul{\bE}\to\ul{\bF}$ from the category
$\abe_{w_0}(X, W)$ we define a mapping cone $\Cone(p)$ as an object
$$
\Cone(p)=\xymatrix{ \Bigl( \bF_1\oplus \bE_0 \ar@<0.6ex>[r]^{c_1} &
\bF_0\oplus \bE_1 \ar@<0.6ex>[l]^{c_0} } \Bigl)
$$
such that
$$
c_0=
\begin{pmatrix}
f_0 & p_1\\
0 & -e_1
\end{pmatrix},
\qquad
c_1=
\begin{pmatrix}
f_1 & p_0\\
0 & -e_0
\end{pmatrix}.
$$
There are  maps
$q: \ul{\bF}\to \Cone(p), \; g=(\id , 0)$ and $r: \Cone(p)\to\ul{\bE}[1],\;
r=(0, -\id).$

The standard triangles in the category
$\dt_{w_0}(X, W)$ are defined to be the triangles of the form
$$
\ul{\bE}\stackrel{p}{\lto} \ul{\bF}\stackrel{q}{\lto} \Cone(p)
\stackrel{r}{\lto} \ul{\bE}[1]
$$
for some $p\in \abe_{w_0}(X, W).$

\begin{definition}
A triangle
$\ul{\bE}{\to} \ul{\bF}{\to} \ul{\bG}
{\to} \ul{\bE}[1]$ in $\dt_{w_0}(X, W)$
is  called  an exact triangle if it is isomorphic to a
standard triangle.
\end{definition}

\begin{proposition}\label{trstr}
The category $\dt_{w_0}(X, W)$ endowed with the translation functor
$[1]$ and the above class of exact triangles becomes a triangulated
category.
\end{proposition}
\begin{dok}  The proof is straighforward.
\end{dok}

We define a triangulated category $\MF_{w_0}(X, W)$ of matrix
factorizations on $(X, W)$ as a Verdier quotient of $\dt(X, W)$ by a
triangulated subcategory of ``acyclic'' objects. This quotient will
also be called a triangulated category of D-branes of type B in the
LG-model $(X, W)$ over $w_0.$

More precisely, for any complex of objects of the category
$\abe_{w_0}(X, W)$
$$
\ul{\bE}^{i}\stackrel{d^i}{\lto}
\ul{\bE}^{i+1}\stackrel{d^{i+1}}{\lto}'
\cdots\stackrel{d^{j-1}}{\lto}\ul{\bE}^{j}
$$
we can consider a totalization $\ul{\bT}$ of this complex. It is a
pair with
\begin{equation}\label{tot}
\bT_1=\bigoplus_{
k+m\equiv 1\mod 2}\bE_k^m,\qquad
\bT_0=\bigoplus_{
k+m\equiv 0 \mod 2}\bE_k^m,\quad k=0,1,
\end{equation}
and with $t_l=d^m_k+(-1)^m e_k$ on the component $\bE_k^m,$ where
$l=(k+m)\mod\; 2.$

Denote by
$\Ac_{w_0}(X, W)$ the minimal full triangulated subcategory that
contains totalizations of all acyclic complexes
in the exact category
$\abe_{w_0}(X,W).$
It is easy to see that $\Ac_{w_0}(X, W)$ coincides with
the minimal full triangulated subcategory containing  totalizations of
all short exact sequences in
$\abe_{w_0}(X,W).$
\begin{definition} We define the triangulated category of matrix
  factorizations $\MF_{w_0}(X, W)$ on $X$ with a
superpotential $W$ as the Verdier quotient
$\dt_{w_0}(X, W)/\Ac_{w_0}(X, W).$
\end{definition}
In particular, this definition implies that any short exact sequence
in $\abe_{w_0}(X, W)$
becomes an exact triangle in $\MF_{w_0}(X, W).$
\begin{remark} {\rm Triangulated categories of this type appeared in
    \cite{Po} under name ``derived categories of the second kind''.
This construction of matrix factorizations is also discussed in
\cite{KKP}.}
\end{remark}

\begin{remark}
{\rm Note that the category of matrix factorization $\MF_{w_0}(X, W)$
defined above has a natural differential graded enhancement obtained
as a DG quotient of the DG category $\dg_{w_0}(X, W)$ by the
corresponding DG subcategory of all totalizations of acyclic
objects. The existence of this DG quotient follows from general
results of Bernhard Keller and Vladimir Drinfeld.}
\end{remark}

\section{Matrix factorizations and categories of singularities}

In this section we discuss a connection between the category of matrix
factorizations $MF_{w_0}(X, W)$ introduced above and the triangulated
category of singularities $\dsing{X_{w_0}}$ of the fiber $X_{w_0}$ of
the map $W$ over a $\kk$\!-point $w_0.$ Without loss of generality we
can take $w_0=0.$ Denote by $X_{0}$ the fiber of $W: X\to \AA^{1}$
over the point $0$ and by $i: X_0\hookrightarrow X$ the closed
embedding.

With any pair $\ul{\bE}$ we can associate a short exact sequence
\begin{equation}\label{shseq}
0\lto \bE_1 \stackrel{e_1}{\lto} \bE_0 \lto \Coker e_1 \lto 0
\end{equation}
of coherent sheaves on $X.$

We can attach to an object $\ul{\bE}$ the sheaf $\Coker e_1.$ This is
a sheaf on $X.$ But the multiplication with $W$ annihilates it. Hence,
we can consider $\Coker e_1$ as a sheaf on $X_{0},$ i.e. there is a
sheaf $\E$ on $X_0$ such that $\Coker e_1\cong i_*\E.$ Any morphism
$p:\ul{\bE}\to\ul{\bF}$ in $\abe_{0}(X, W)$ gives a morphism between
cokernels.  In this way we get a functor $\Cok:\abe_{0}(X,W)\lto
\coh(X_{0}).$

\begin{lemma}\label{isom} Let $\ul{\bE}$ and $\ul{\bF}$ be two pairs
  on $X.$ Then we have
\begin{itemize}
\item[a)] If $p:\ul{\bF}\to \ul{\bE}$ is a morphism which
  induces an isomorphism between  the sheaves $\Coker f_1$ and
$\Coker e_1,$ then $p$ becomes an isomorphism in $\MF_0(X, W).$
\item[b)] For any morphism $a: \Coker f_1\to \Coker e_1$ there is a
  pair $\ul{\bF}'$ together with  morphisms $p:\ul{\bF}'\to \ul{\bE}$
  and
$s:\ul{\bF}'\to \ul{\bF}$ such that $\Cok(s)$ is an isomorphism and
  $a=\Cok(p)\Cok(s)^{-1}.$ Moreover, if $a$ is a surjection, then
$\ul{\bF}'$ can be chosen so that $p$ is a surjection as well.
\item[c)] If the sheaves $\Coker f_1$ and
$\Coker e_1$ are isomorphic, then $\ul{\bF}$ and $\ul{\bE}$ are
  isomorphic in
$\MF_0(X, W).$
\end{itemize}
\end{lemma}
\begin{dok}
a) Let us take a locally free covering $\bG\to\bE_1$ and consider a
trivial pair $\ul{\bG}$ with $\bG_1=\bG_0=\bG$ and $g_1=\id, g_0=W.$
We have a map of pairs $\ul{\bG}\to\ul{\bE}$ that induces a natural
map $p': \ul{\bF}\oplus\ul{\bG}\to\ul{\bE}.$ If $p:\ul{\bF}\to
\ul{\bE}$ gives an isomorphism between sheaves $\Coker f_1$ and
$\Coker e_1$ then $p'$ is a surjection and the kernel of $p'$ is a
$0$\!-homotopic pair. Hence, $p'$ and $p$ are isomorphisms in
$\MF_0(X, W).$

b) A morphism $a:\Coker f_1\to \Coker e_1$ induces a morphism
$a':\bF_0\to\Coker e_1.$ Consider a pull back
of $a'$ along the surjection $\bE_0\to \Coker e_1$ and denote by $\bU$
a Cartesian product

$$
\xymatrix{
\bU\ \ar[d]\ar[r]&\bF_0\ar[d]^{a'}\\
\bE_0\ar[r]& \Coker e_1
}
$$ Since $X$ has enough locally free sheaves we can find a locally
free sheaf $\bF_0'$ with a surjection on $\bU.$ It induces a map $p_0:
\bF_0'\to\bE_0$ and a surjection $s_0:\bF_0'\twoheadrightarrow \bF_0.$

The surjection $s_0$ gives a map $\bF_0'\to \Coker f_1.$ Denote by
$f_1': \bF_1'{\to}\bF_0'$ the kernel of the latter map.  The
surjection $s_0$ defines a surjection $s_1: \bF_1'\to \bF_1$ and the
sheaf $\bF_1'$ is locally free as extension of $\bF_1$ with the kernel
of $s_0.$ Since the multiplication by $W$ acts trivially on $\Coker
f_1$ there is a map $f_0': \bF_0'\to \bF_1'$ such that the composition
$f_1'f_0'$ coincides with the multiplication by $W$ on $\bF_0'.$ It is
easy to see that the composition $f_0' f_1'$ is equal to
multiplication with $W$ on $\bF_1'$ as well.  Thus, we obtain a pair
$\ul{\bF}'$ and a surjection $s: \ul{\bF}'\twoheadrightarrow \ul{\bF}$
for which $\Cok(s)$ is an isomorphism.

On the other hand, we have a map $p_0: \bF_0'\to\bE_0$ that can be
included in the following commutative square
$$
\xymatrix{
\bF_0'\ \ar[d]_{p_0}\ar[r]&\Coker f_1\ar[d]^{a}\\
\bE_0\ar[r]& \Coker e_1
}
$$ This commutative diagram can be extended (uniquely) to a map of the
pairs $p:\ul{\bF}'\to\ul{\bE}.$ If now $a$ was a surjection, then
$p_0:\bF_0'\to\bE_0$ is also surjection. As in a) above take a locally
free covering $\bG\to\bE_1$ and consider a trivial pair $\ul{\bG}$
with $\bG_1=\bG_0=\bG$ and $g_1=\id, g_0=W.$ We have a map of pairs
$\ul{\bG}\to\ul{\bE}.$ Changing $\ul{\bF}'$ to the direct sum
$\ul{\bF}'\oplus\ul{\bG}$ we obtain a surjection on $\ul{\bE}.$ Thus,
b) is proved.  The statements b) and a) imply c) directly.
\end{dok}

One of the main benefits of this lemma is that it allows us to replace
a given pair $\ul{\bE}$ by pairs $\ul{\bE}'$ which are isomorphic to
$\ul{\bE}$ in the category $\MF_{0}(X, W)$ but are better suited to
our purposes.  Indeed, consider a pair $\ul{\bE}$ and the coherent
sheaf $\E\in \coh(X_0)$ that is the cokernel of the map $e_1.$ Since
$X$ has enough locally free sheaves we can find a locally free sheaf
$\bL$ on $X$ with a surjection $i^* \bL\twoheadrightarrow \E.$ By
Lemma \ref{isom} b) there is a pair $\ul{\bF}$ such that $\Coker
f_1\cong i_*i^*\bL$ and a surjection
$p:\ul{\bL}\twoheadrightarrow\ul{\bE}$ that induces the surjection
$i^* \bL\twoheadrightarrow \E.$ Taking the kernel of $p$ we obtain a
pair $\ul{\bG}$ that is isomorphic to $\ul{\bE}[1]$ in $\MF_{0}(X,
W).$ In particular it follows that $\ul{\bL}$ is the
$0$\!-object in $\MF_{0}(X, W)$ and that the sequence
$\ul{\bG}\to\ul{\bL}\to\ul{\bE}$ gives an exact triangle in
$\MF_{0}(X, W).$ Repeating this procedure for $\ul{\bG}$ we get
$\ul{\bE}'\cong \ul{\bG}[1]\cong\ul{\bE}[2]\cong\ul{\bE}.$

With these preliminaries we are ready to show that the functor $\Cok$
induces exact functors from the triangulated categories $\dt_{0}(X,
W)$ and $\MF_0(X, W)$ to the triangulated category of singularities
$\dsing{X_{0}}$ of the fiber $X_0.$

\begin{proposition} The functor $\Cok: \abe_{0}(X, W) \to \coh(X_{0})$
induces  exact functors \linebreak
$\Pi: \dt_{0}(X, W)\to \dsing{X_{0}}$ and $\Sigma:
 \MF_{0}(X, W)\to \dsing{X_{0}}$ between triangulated categories.
\end{proposition}
\begin{dok}
We have a functor $\abe_{0}(X, W)\to \dsing{X_0}$ which is the
composition of $\Cok$ and the natural functor from $\coh(X_0)$ to
$\dsing{X_0}.$ To prove the existence of the functor $\Pi$ we must
show that any morphism $p=(p_1, p_0): \ul{\bE} \to \ul{\bF}$ that is
homotopic to $0$ goes to the $0$\!-morphism in $\dsing{X_0}.$ Fix a
homotopy $(s_1, s_0),$ where $s_1:\bE_1\to \bF_0$ and $s_0: \bE_0\to
\bF_1.$ Consider the following decomposition of $p$

$$
\xymatrix{
\bE_1 \ar[d]_{\binom{s_1}{p_1}} \ar[r]^{e_1}  & \bE_0
 \ar[d]^{\binom{s_0}{p_0}} \ar[r]& i_*\E \ar[d]
\\
\bF_0\oplus \bF_1 \ar[d]_{pr} \ar[r]^{c_1} & \bF_1\oplus
\bF_0\ar[d]^{pr}\ar[r] & i_* i^* \bF_0 \ar[d]
 & {}\save[]*\txt{where}\restore&
{c_1=
\begin{pmatrix}
-f_0 & \id_{\bF_1}\\
0 & f_1
\end{pmatrix}.
}
\\
\bF_1 \ar[r]^{f_1} & \bF_0 \ar[r] &
i_*\F
}
$$ This gives the decomposition of $\Cok(p)$ through a locally free
object $i^* \bF_0.$ Hence, $\Cok(p)=0$ in the category $\dsing{X_0}$
and we obtain a functor $\Pi$ from $\dt_{0}(X, W)$ to $\dsing{X_{0}}.$
It is easy to check that the functor $\Pi$ is exact.

To see that the functor $\Pi$ induces a functor from the quotient
category $\MF_{0}(X, W)$ to $\dsing{X_{0}}$ we should to check that
for any object $\ul{\bT}\in \Ac_0(X, W)$ the cokernel of $t_1$ is
locally free on $X_0.$ The subcategory $\Ac_0(X, W)$ is generated by
totalizations of short exact sequences in $\abe_{0}(X,W).$ Let us
consider a short exact sequence
\begin{equation}\label{short}
0\lto\ul{\bG}\stackrel{q}{\lto}\ul{\bE}\stackrel{p}{\lto}\ul{\bF}\lto 0
\end{equation}
and denote by $\ul{\bT}$ its totalization as in (\ref{tot}). We also
consider a complex
$$
\begin{CD}
0 @>>> \bG_0 @>{\left(\begin{smallmatrix}q_0\\g_0\end{smallmatrix}\right)}>>
\bE_0\oplus\bG_1 @>{\left(\begin{smallmatrix}p_0 & 0\\- e_0 &
    q_1\end{smallmatrix}\right)}>> \bF_0\oplus \bE_1 @>(f_0,
p_1)>>\bF_1 @>>> 0
\end{CD}
$$ that is a totalization of the complex (\ref{short}) in the category
of complexes of coherent sheaves on $X.$ This complex is acyclic.  We
denote by $\bU$ the image of the middle map and by $\psi=(\psi_0,
\psi_1): \bE_0\oplus\bG_1\to\bU$ and $\phi:\bU\to\bF_0\oplus \bE_1 $
the canonical surjection and injection. We have the following
commutative diagram

$$
\xymatrix{
\bE_0\oplus \bG_1 \ar[d]_{can}\ar[r]^{\left(\begin{smallmatrix}\psi_0
    & \psi_1\\0 & g_1\end{smallmatrix}\right)} &\bU\oplus \bG_0
\ar[d]^{\phi\oplus\id}\ar[r] &
i_* \H \ar[d]
\\
\bT_1=\bF_1\oplus\bE_0\oplus \bG_1 \ar[d]_{pr} \ar[r]^{t_1} &
\bT_0=\bF_0\oplus\bE_1\oplus
\bG_0\ar[d]^{(f_0, p_1, 0 )}\ar[r] & i_*\T \ar[d]
 & {}\save[]*\txt{where}\restore&
{t_1=
\begin{pmatrix}
f_1 & p_0 & 0\\
0 & -e_0 & q_1 \\
0& 0 & g_1
\end{pmatrix}
}
\\
\bF_1  \ar[r]^{W}  & \bF_1
  \ar[r]& i_* i^* \bF_1
}
$$ in which all columns and all rows are short exact sequences on $X.$
Hence, the object $\T$ is an extansion of $i^*\bF_0$ and $\H.$ On the
other hand, we also have the commutative diagram

$$
\xymatrix{
\bG_0\
\ar[d]_{\left(\begin{smallmatrix}q_0\\
g_0\end{smallmatrix}\right)}\ar[r]^{W}&\bG_0\ar[d]^{can}\ar[r]
&
i_* i^*\bG_0 \ar@{=}[d]\\
\bE_0\oplus \bG_1 \ar[r]^{\left(\begin{smallmatrix} \psi_0 & \psi_1\\
    0 & g_1\end{smallmatrix}\right)}\ar[d]_{\psi} &\bU\oplus \bG_0
\ar[d]^{pr}\ar[r] &
i_* \H
\\
\bU\ar@{=}[r]& \bU&
}
$$
which shows that $\H\cong i^*\bG_0$ is locally free on $X_0.$
Thus, the object $\T,$ which is the cokernel of the map $t_1,$ is an extension of two locally free sheaves
$i^*\bF_1$ and $i^*\bG_0$ on $X_0.$ Therefore, it is locally free as well.
\end{dok}

\begin{proposition}\label{embob}
If $\Sigma(\ul{\bE})\cong 0,$ then $\ul{\bE}\cong0$ in the category
$\MF_{0}(X, W).$
\end{proposition}
\begin{dok}
Let  $\E=\Sigma(\ul{\bE}).$ If $\E\cong 0$ in $\dsing{X_0}$ then $\E$
is a perfect complex.
On the other hand, we have 2-periodic resolution of the form
$$
\{\cdots\lto i^* \bE_0 \lto i^*\bE_1 \lto i^*\bE_0\}\lto\E\lto 0
$$
Since $\E$ is perfect then the kernel of the first map of the brutal
truncation $\sigma^{\ge 2n+1}$ of this resolution for sufficient
negative $n$ is a locally free sheaf. On the other hand it is isomorphic
to $\E.$ Hence, $\E$ is locally free if it is perfect.

Note that if $\E\cong i^*\bF,$ then by Lemma \ref{isom} the object
$\ul{\bE}$ is isomorphic in the category $\MF_{0}(X, W)$ to the object
$\ul{\bF}$ with $\bF_1\cong\bF_0\cong\bF$ and $f_1=W, f_0=\id.$ This
implies that $\ul{\bE}\cong 0$ in the category $\MF_{0}(X, W),$
because $\ul{\bF}$ is $0$\!-homotopic.

If now a locally free sheaf $\E$ on $X_0$ has a right resolution of
the form
\begin{equation}\label{rightres}
0\lto\E\lto i^* \bG^{-k} \lto i^*\bG^{-k+1} \lto \cdots \lto
i^*\bG^0\lto 0,
\end{equation}
then by Lemma \ref{isom} b) we can construct an exact sequence of
pairs of the form
\begin{equation}\label{respai}
0\lto \ul{\bE}'\lto
\ul{\bL}^{-k}\lto\ul{\bL}^{-k+1}\cdots\lto\ul{\bL}^0\lto 0
\end{equation}
such that the functor $\Cok$ sends the sequence (\ref{respai}) to the
 sequence (\ref{rightres}).
 As we noted above, all pairs
$\ul{\bL}^{-i}$ are isomorphic to $0$ in $\MF_0(X, W).$ Hence, the
 object $\ul{\bE}'\cong 0$ in $\MF_0(X, W)$ too.
And by Lemma \ref{isom} c) the object $\ul{\bE}$ is isomorphic to $0$
 in the category  $\MF_0(X, W)$ as well.

Finally, consider a general case. Let $\E=\Sigma(\ul{\bE})$ be a
locally free sheaf on $X_0.$
Let us consider the following complex
\begin{equation}\label{perres}
\underset{2m}{\underbrace{i^*\bE_1\to\cdots \to i^*\bE_1 \to i^*\bE_0}}
\end{equation}
which is concentrated in degrees $(-2m+1)$ to $0$
for $2m>\dim X_0.$ It has two nontrivial cohomology $H^{-2m+1}$ and
$H^0,$ both of which are isomorphic  to $\E.$
Since $\E$ is locally free and $2m>\dim X_0$ this complex in $\db{\coh X_0}$
is isomorphic to the direct sum of its cohomology. Hence,  there is a
map from $\E$ to the complex (\ref{perres}) in $\db{\coh X_0}.$
Therefore, we can find a locally free resolution $P^{\cdot}$ of $\E$
of the form
\begin{equation}\label{resol}
\{\cdots\lto i^* \bP^{-k}\lto\cdots\lto i^* \bP^0\}\lto \E
\end{equation}
and a map from this resolution to the complex (\ref{perres}) that acts
identically on the 0-th cohomology.
Take the brutal truncation $\sigma^{\ge -2m+2}(P^{\cdot})$ of the
resolution (\ref{resol}) and consider the corresponding composition
map from the complex $\sigma^{\ge -2m+2}(P^{\cdot})$ to
(\ref{perres})
$$
\xymatrix{
&0\ar[r]\ar[d]&i^*\bP^{-2m+2} \ar[d]_{} \ar[r]  &\cdots \ar[r] &
  i^*\bP^0\ar[r]\ar[d]&0
 \ar[d]^{}
\\
0\ar[r]&i^*\bE_1\ar[r]& i^*\bE_0\ar[r]&\cdots \ar[r] & i_* \bE_0
\ar[r]&0
}
$$
Since this map induces an isomorphism on the 0-th cohomology,  the cone of
this map
\begin{equation}\label{totali}
i^*\bE_1\oplus i^*\bP^{-2m+2}\lto\cdots\lto i^*\bE_1\oplus
i^*\bP^0\lto i^* \bE_0
\end{equation}
has only one nontrivial cohomology $H^{-2m+1}.$ Denote this cohomology
by $\T.$ It is easy to see that $\T$ is isomorphic to the direct
sum
of $\E$ and $H^{-2m+2}(\sigma^{\ge -2m+2} P^{\cdot}).$ This follows
from the fact that in the derived category
$\db{\coh X_0}$ the complex (\ref{totali}) is a cone of the map \linebreak
$\sigma^{\ge -2m+2}(P^{\cdot})\to \E\oplus\E[2m-1]$ which factors through
$\E.$

Thus, the locally free sheaf $\E$ is a direct summand of the locally
free sheaf $\T$ that has a right resolution of the form
(\ref{totali}).  It was shown above that any pair $\ul{\bT},$ for
which $\Cok(\ul{\bT})\cong \T,$ is isomorphic to $0$ in the category
$\MF_0(X, W).$ Now consider the maps
$\E\stackrel{i}{\to}\T\stackrel{\pi}{\to}\E,$  whose composition is
the identity. By Lemma \ref{isom} b) there is a pair $\ul{\bT}$ and a
surjection $p:\ul{\bT}\to\ul{\bE}$ such that $\Cok(p)=\pi.$ Applying
again Lemma \ref{isom} b) we find a pair $\ul{\bE}'$ and a map
$q:\ul{\bE}'\to \ul{\bT}$ such that $\Cok(q)=i.$ This implies that the
pairs $\ul{\bE}'$ and $\ul{\bE}$ are isomorphic in $\MF_0(X, W)$ and
are isomorphic to a direct summand of the object $\ul{\bT}.$ Since
$\ul{\bT}$ is isomorphic to the 0-object we have that $\ul{\bE}\cong
0$ in the category $\MF_0(X, W)$ as well.
\end{dok}

\begin{theorem}\label{main2} Let $X$ be a scheme that satisfies
  condition (ELF). Then
the natural functor $\Sigma: \MF_{0}(X, W)\to \dsing{X_{0}}$ is fully
faithful.
\end{theorem}
\begin{dok}
First we will show that $\Sigma$ is full.  Let $\F$ and $\E$ be
sheaves on $X_0$ that are cokernels of $f_1$ and $e_1$ for two pairs
$\ul{\bF}$ and $\ul{\bE},$ respectively.  By the definition of a
localization any morphism $u$ from $\F$ to $\E$ in $\dsing{X_0}$ can
be represented by a pair of morphisms in $\db{\coh(X_0)}$ of the form
\begin{equation}\label{domik}
\F\stackrel{s}{\longleftarrow}A\stackrel{a}{\lto} \E
\end{equation}
such that the cone $C(s)$ is a perfect complex.  The two cohomology
sheaves of the complex $i^*\bE_1\to i^*\bE_0$ are isomorphic to $\E$
and we obtain a canonical map $\phi:\E\to \E[2]$ in $\db{\coh(X)}$
that becomes an isomorphism in the quotient category $\dsing{X_0}.$
The morphism $u$ induces a morphism $u_k$ from $\F$ to $\E[2k]$ in
$\dsing{X_0}$ that is represented by the roof
$$
\F\stackrel{s}{\longleftarrow}A\stackrel{a_k}{\lto} \E[2k]
$$
with $a_k=\phi^k a.$ Since the cone of $s$ is perfect, it is
isomorphic to a bounded complex of locally free sheaves.  This means
that for a sufficiently large $k$ there is no nontrivial map from $C(s)$
and $C(s)[-1]$ to $\E[2k].$ Hence, the map $u_k$ in $\dsing{X_0}$ is
represented by a map $\tilde{u}$ from $\F$ to $\E[2k]$ in
$\db{\coh(X_0)}.$

Now we take a surjective morphism $\pi: i^* \bG^0\twoheadrightarrow\F$
that erases the map $\tilde{u}:\F\to \E[2k].$ Denote by $\H$ the
kernel of $\pi.$ The map $\tilde{u}$ induces a map
$\tilde{\tilde{u}}:\H\to\E[2k-1].$ By the same argument we find a
surjective morphism $i^* \bG^{-1}\twoheadrightarrow \H$ that erases
the map $\tilde{\tilde{u}}:\H\to \E[2k-1].$

Repeating the above procedure we construct an acyclic complex of the form
\begin{equation}\label{Lres}
0\lto\F'\lto i^* \bG^{-2k+1}\lto\cdots\lto i^* \bG^0\lto \F\lto 0
\end{equation}
and a map $u':\F'\to\E$ such that $\tilde{u}$ is the composition of
the canonical map $\F\to\F'[2k]$ and $u'[2k]:\F'[2k]\to\E[2k].$

By Lemma \ref{isom} b) we can construct an exact sequence of pairs of the form
\begin{equation}\label{Lpair}
0\lto \ul{\bF}'\lto \ul{\bL}^{-2k+1}\lto\cdots\lto\ul{\bL}^0\lto
\ul{\bF}\lto 0
\end{equation}
such that the application of functor $\Cok$ to the sequence
(\ref{Lpair}) gives the sequence (\ref{Lres}). All pairs
$\ul{\bL}^{-i}$ are isomorphic to $0$ in $\MF_0(X, W).$ Hence, the
object $\ul{\bF}'\cong \ul{\bF}$ in $\MF_0(X, W).$
Moreover, there is an equality $u'\cdot \Sigma(\alpha)=u$ of morphisms
in $\dsing{X_0},$ where
 $\alpha:\ul{\bF}\stackrel{\sim}{\to}\ul{\bF}'$ is the isomorphism in
the category $\MF_0(X, W).$

 Thus, it is enough to show that the morphism $u':\F'\to\E$ is equal
 to $\Sigma(q)$ for some map $\beta: \ul{\bF}'\to \ul{\bE}$ in the
 category $\MF_0(X, W).$ By Lemma \ref{isom} b) there is a pair
 $\ul{\bF}''$ and morphisms of pairs $p:\ul{\bF}''\to\ul{\bE}$ and $s:
 \ul{\bF}''\to\ul{\bF}'$ such that $\Cok(s)$ is an isomorphism and
 $u'=\Cok(p)\cdot\Cok(s)^{-1}.$ Again by Lemma \ref{isom} a) the morphism
 $s$ becomes an isomorphism in the category $\MF_0(X, W).$ Hence, for
 $\beta=p s^{-1}$ in $\MF_0(X, W)$ we obtain that $\Sigma(\beta)=u'.$
 Therefore, the functor $\Sigma$ is full.

Now we prove that $\Sigma$ is faithful. It is a standard statement
asserting that a full exact functor between triangulated categories
with trivial kernel is also faithful.  Indeed, let $p: \ul{\bE}\to
\ul{\bF}$ be a morphism for which $\Sigma(p)=0.$ Complete $p$ to an
exact triangle
$$
\ul{\bE}\stackrel{p}{\lto}\ul{\bF}\stackrel{q}{\lto}\ul{\bG}\lto
\ul{\bE}[1].
$$
Then the identity map of $\Sigma(\ul{\bF})$ factors through the map
$\Sigma(\ul{\bF})\stackrel{\Sigma(q)}{\lto} \Sigma(\ul{\bG}).$ Since
$\Sigma$ is full, there is a map $s:\ul{\bF}\to \ul{\bF}$ factoring
through $q: \ul{\bF}\to \ul{\bG}$ such that $\Sigma(s)=\id.$ Hence,
the cone $C(s)$ of the map $s$ goes to zero under the functor $\Sigma.$ By
Proposition \ref{embob} the object $C(s)$ is the zero object too,
so $s$ is an isomorphism.  Thus, $q:\ul{\bF}\to\ul{\bG}$ is a split
monomorphism and $p=0.$
\end{dok}

In general the functor $\Sigma$ is not necessarily an
equivalence and it is very interesting to understand a difference
between $\MF_0(X, W)$ and $\dsing{X_0}$ for singular LG-models.
However, for regular schemes we have an equivalence.

\begin{theorem} Let $X$ be a scheme that satisfies condition (ELF). If
  X is regular then the functor
$\Sigma: \MF_{0}(X, W)\to \dsing{X_0}$ is an equivalence of
  triangulated categories.
\end{theorem}
\begin{dok}
As we showed above the functor $\Sigma$ is fully faithful. To
complete the proof that $\Sigma$ is an equivalence we need to check
that every object $A\in \dsing{X_0}$ is isomorphic to
$\Sigma(\ul{\bE})$ for some pair $\ul{\bE}.$ By Proposition
\ref{fint}, since $X_0$ is Gorenstein  any object $A\in \dsing{X_0}$
is isomorphic to the image of a coherent sheaf $\E$ such that
$\ext^{i}_{X_0}(\E, \O_{X_0})=0$ for all $i>0.$ Consider an
epimorphism $\bE_0\twoheadrightarrow i_* \E$ of sheaves on $X$  with
locally free $\bE_0.$ Denote by $e_1: \bE_1\to \bE_0$ the kernel of
this map.
 Since the multiplication with $W$ gives the zero map on $\E,$
there is a map $e_0: \bE_0\to \bE_1$ such that
$e_0 e_1=W$ and $e_1 e_0=W.$ We get a pair
$$
\ul{\bE}:=
\Bigl(
\xymatrix{
\bE_1 \ar@<0.6ex>[r]^{e_1} &\bE_0 \ar@<0.6ex>[l]^{e_0}
}
\Bigl)
$$
and we need only to check that $\bE_1$ is locally free. It  follows
from the fact that for any closed point $t:x\hookrightarrow X$ we
have
\begin{equation}\label{last}
\Ext^i_X(\bE_1, t_* \O_x)= 0
\end{equation}
for all $i>0.$ To show it we note that by Lemma \ref{rres}  the
sheaf $\E$ has a right locally free resolution on $X_0.$ For any
local free sheaf $\P$ on $X_0$ we have $ \Ext^i_X(i_*\P, t_* \O_x)=
0 $ for $i>1.$ Since $X$ is regular, the abelian category of coherent
sheaves on $X$ has finite cohomological dimension. Therefore,  we
obtain    $\Ext^i_X(i_*\E, t_* \O_x)= 0$ for $i>1.$ This
implies (\ref{last}) and the theorem.
\end{dok}

\end{document}